\documentclass[manmat,francais,envcountsame,envcountreset,numbook]{svjour}
\usepackage{amssymb,amsmath}
\renewcommand{\refname}{R\'ef\'erences}
\spnewtheorem{rem}[theorem]{Remarque}{\bf}{\it}
\newcommand{\hfl}[2]{\smash{\mathop{\hbox to 7
mm{\rightarrowfill}}\limits^{\scriptstyle #1}_{\scriptstyle #2}}}
\newcommand{\vfl}[2]{\llap{$\scriptstyle #1$}\left\downarrow\vbox to
6mm{}\right.\rlap{$\scriptstyle#2$}}
\def\limproj{\mathop{\oalign{$\mathrm{lim}$\cr\hidewidth$\longleftarrow$\hidewidth\cr}}}
\def\limindu{\mathop{\oalign{$\mathrm{lim}$\cr\hidewidth$\longrightarrow$\hidewidth\cr}}}
\title{Sur la torsion de la distribution ordinaire universelle attach\'ee \`a un corps de
nombres}
\date{}
\titlerunning{torsion de la distribution universelle}
\author{Jean-Robert BELLIARD\inst{1}\and
 Hassan OUKHABA\inst{2}}
\authorrunning{J.-R. Belliard, H. Oukhaba}
\institute{1 School of Mathematical Sciences,
University of Nottingham\\
2 UMR 6623,Universit\'e de Franche-Comt\'e, Besan\c 
con\\
  \email{pmxjrb@maths.nott.ac.uk}
\hfill \email {hassan@math.univ-fcomte.fr}}

\begin{document}
\maketitle
\setcounter{section}{-1}
\bibliographystyle{plain}

\begin{abstract} We study the torsion subgroup of the universal
ordinary distribution related to a general number field. We describe a
way to control this subgroup. We apply this method to the special case
of an imaginary quadratic field, and we give examples of such fields
where these torsion subgroups are non-trivial.
\end{abstract}

\section{Introduction}\label{intro}
Kubert dans son article \cite {Ku79} a introduit le formalisme 
des distributions sur $\mathbb{Q}$. Il y d\'emontre aussi 
que la distribution ordinaire universelle sur $\mathbb{Q}$
 est sans $\mathbb{Z}$-torsion.
La question de l'existence de torsion dans 
la distribution ordinaire universelle $\mathcal{A}_K$ attach\'ee \`a un corps 
de nombres $K$ quelconque se pose.
Si on \'ecrit $\mathcal{A}_K$ comme limite inductive de
ses sous groupes
de niveau $\mathcal{A}_{\mathfrak n}$ index\'es
par les id\'eaux entiers de $K$, alors Yin 
 a montr\'e que
$\mathcal{A}_{\mathfrak n}$ est sans $\ell$-torsion pourvu que le nombre premier $\ell$ ne divise
pas l'ordre du
groupe de Galois $G_{\mathfrak n}$ (voir \cite{Y00}).

Dans
cet article nous
nous int\'eressons au cas o\`u $\ell$ divise cet ordre et nous d\'egageons 
des conditions suffisantes assez fines pour que $\mathcal{A}_{\mathfrak n}$
 soit sans $\ell$-torsion. Ces conditions sont li\'ees \`a la
structure du $\ell$-Sylow d'un sous-groupe de $G_{\mathfrak n}$.
 Avec comme objectif les applications arithm\'etiques de \cite{Ou98} 
nous appliquons ensuite cette d\'emarche au cas o\`u le
corps de base $K$ est un corps quadratique imaginaire.
 Nous obtenons ainsi une majoration pour l'ordre
de la torsion des sous-groupes de niveau.
Nous montrons ensuite que la distribution ordinaire universelle 
attach\'ee a tout corps quadratique
non-ramifi\'e en $2$ contient de la torsion
(en fait dans ces cas la borne calcul\'ee pr\'ec\'edemment est
atteinte d\`es que le corps quadratique est principal).

\section{Formalisme g\'en\'eral}\label{nota}

Avant d'entamer notre \'etude nous donnons ci-dessous quelques-unes des
 notations utilis\'ees dans la suite.

\begin{list}{}{}
\item $\mathcal{O}_K$:= l'anneau des entiers de $K$.
\item $\mathfrak{m}_\infty$:= le produit des places r\'eelles de $K$.
\item $\Psi_0$:= Le mono\"{\i}de des id\'eaux entiers non nuls de $K$.
\item $\bar \Psi_0$:= Le groupe des id\'eaux fractionnaires de $K$.
\end{list}
De plus pour tout $\mathfrak{n}\in \Psi_0$ nous noterons :
\begin{list}{}{} 
\item $I_{\mathfrak n}$:= le groupe des id\'eaux fractionnaires de K premiers avec
$\mathfrak{n}$.
\item $\mathcal{R}_{\mathfrak n}$:= le sous groupe de $I_{\mathfrak n}$ form\'e des id\'eaux
principaux
$x\mathcal{O}_K$ tels que $x\equiv1$ $\mathrm{mod}^\times\mathfrak{n}$ et $\sigma(x)>0$ pour tout
plongement r\'eel
$\sigma$ de $K$.
\item $K_{\mathfrak n}$:=le corps de rayon modulo $\mathfrak{nm}_\infty$.
\item $T_{\mathfrak p}(\mathfrak{n})$:= le groupe d'inertie de $\mathfrak{p}$ dans $K_{\mathfrak
n}/K$,
$\mathfrak{p}$ \'etant un id\'eal premier.
\item $G_{\mathfrak n}:=\mathrm{Gal}(K_{\mathfrak n}/K)$.
\end{list}\par
Rappelons que le groupe $G_{\mathfrak n}$ est isomorphe \`a $I_{\mathfrak
n}/\mathcal{R}_{\mathfrak n}$, par l'isomorphisme de r\'eciprocit\'e d'Artin.
On notera $(\mathfrak{a},K_{\mathfrak n}/K)$
l'automorphisme de
$K_{\mathfrak n}/K$ associ\'e \`a l'id\'eal $\mathfrak{a}\in I_{\mathfrak n}$ par cet
isomorphisme. Si $X$ est une
partie de $G_{\mathfrak n}$ alors on posera $s(X):=\sum\sigma\in\mathbb{Z}[G_{\mathfrak n}]$,
o\`u $\sigma$
parcourt les \'el\'ements de $X$ et $\# X$ d\'esignera le cardinal de $X$. Consid\'erons
maintenant les relations d'\'equivalence $\sim_{\mathfrak{n}}$ d\'efinies
sur $\Psi_0$ comme suit :
$$\mathfrak{a}\ \sim_{\mathfrak{n}}\ \mathfrak{b}\iff
pgcd(\mathfrak{a},\mathfrak{n})=pgcd(\mathfrak{b},\mathfrak{n})\quad\mathrm{et}\quad\mathfrak{b}\mathfrak{a}^{-1}\in\mathcal{R}_{\mathfrak{n}^\prime}$$
o\`u $\mathfrak{n}^\prime$ est tel que $\mathfrak{n}^\prime
pgcd(\mathfrak{a},\mathfrak{n})=\mathfrak{n}$. Lorsque
les id\'eaux $\mathfrak{a}$ et $\mathfrak{b}$ sont premiers avec $\mathfrak{n}$ alors
$\mathfrak{b}\sim_{\mathfrak
n}\mathfrak{a}$ si, et seulement si, $(\mathfrak{b},K_{\mathfrak n}/K)=(\mathfrak{a},K_{\mathfrak
n}/K)$. Notons $\psi_{\mathfrak n}(\mathfrak{a})$ la classe de $\mathfrak{a}$ modulo la relation
$\sim_{\mathfrak{n}}$ et $\Psi_{\mathfrak n}:= \Psi_0/\sim_{\mathfrak{n}}$ l'ensemble des classes
d'\'equivalence de cette relation. Il est clair que $\Psi_{\mathfrak n}$ est fini. C'est aussi un
mono\"{\i}de unitaire commutatif pour la multiplication des id\'eaux. 
L'ensemble des \'el\'ements inversibles de $\Psi_{\mathfrak n}$ forme
ainsi un groupe, qui est naturellement isomorphe \`a $G_{\mathfrak n}$.
Si $\mathfrak{n}\mid\mathfrak{m}$ alors
on a une surjection $\Psi_{\mathfrak{n},\mathfrak{m}}\colon\Psi_{\mathfrak
m}\longrightarrow\Psi_{\mathfrak n}$. La relation $\sim_{\mathfrak n}$ a \'et\'e  introduite par
Deligne et Ribet dans leur article \cite{DR80} afin
d'\'etudier les valeurs des fonctions $L$ des corps de nombres
totalement r\'eels aux entiers n\'egatifs.\par
L'ensemble $\Psi:=\limproj\Psi_{\mathfrak n}$ est un espace topologique compact totalement
discontinu. Par d\'efinition, cf. \cite{MSD}, une distribution sur $\Psi$ \`a valeurs dans un groupe
ab\'elien $A$ est une fonction additive
$$\{\mathrm{ouverts\ compacts\ de\ }\Psi\}\hfl{\mu}{}A.$$
La projection $g_{\mathfrak n}\colon\Psi\longrightarrow\Psi_{\mathfrak n}$ est continue et surjective.
L'image r\'eciproque $g^{-1}_{\mathfrak n}(\psi_{\mathfrak n}(\mathfrak{a}))$ est donc un ouvert
compact pour tout id\'eal entier $\mathfrak{a}\in\Psi_0$. Posons alors $\mu_{\mathfrak
n}(\mathfrak{a}):=\mu(g^{-1}_{\mathfrak n}(\psi_{\mathfrak n}(\mathfrak{a})))$. On dira que $\mu$
est une distribution ordinaire si $\mathfrak{a}\mathfrak{n}^{-1}=\mathfrak{b}\mathfrak{m}^{-1}$
entra\^{\i}ne $\mu_{\mathfrak n}(\mathfrak{a})=\mu_{\mathfrak m}(\mathfrak{b})$. Dans ce cas on
peut voir $\mu$ comme une application $\bar{\Psi}_0\longrightarrow A$ telle que l'on a 
\begin{equation}\label{mesure}
\mu(\mathfrak{au}^{-1})=\sum_{\mathfrak{b}\in
w(\mathfrak{u},\mathfrak{v},\mathfrak{a})}\mu(\mathfrak{bv}^{-1})
\end{equation}
pour tous $\mathfrak{a},\mathfrak{u},\mathfrak{v}\in \Psi_0$ tels que
$\mathfrak{u}\mid\mathfrak{v}$, o\`u
$w(\mathfrak{u},\mathfrak{v},\mathfrak{a})$ est n'importe quel syst\`eme de repr\'esentants dans $\Psi_0$ de l'image r\'eciproque 
$\Psi_{\mathfrak{u},\mathfrak{v}}^{-1}(\psi_{\mathfrak u}(\mathfrak{a}))$.
 En particulier si on pose $\mathfrak{n}^\prime
pgcd(\mathfrak{a},\mathfrak{n})=\mathfrak{n}$ et $\mathfrak{a}^\prime
pgcd(\mathfrak{a},\mathfrak{n})=\mathfrak{a}$, alors 
$\mu(\mathfrak{an}^{-1})$ ne
d\'epend que de la classe de $\mathfrak{a}^\prime$ dans
$I_{\mathfrak{n}^\prime}/\mathcal{R}_{\mathfrak{n}^\prime}$.\par
La notion de distribution 
intervient fr\'equemment en arithm\'etique.
On peut par exemple trouver diverses r\'ealisations dans 
les articles \cite{DR80}, \cite{Ku79}, \cite{MSD},
 \cite {Ou98}, \cite{Wa}, \cite {Y00}, etc ...  

\section{Distribution ordinaire universelle sur $K$}\label{thmgen}
Notons $[\bar \Psi_0]$ le groupe ab\'elien libre sur $\mathbb{Z}$ engendr\'e par les symboles
$[\mathfrak{a}],\mathfrak{a}\in \bar \Psi_0$. Soit alors $\mathcal{A}$ le groupe $[\bar \Psi_0]$
modulo les relations
\begin{equation}\label{dist}
[\mathfrak{au}^{-1}]-\sum_{\mathfrak{b}\in
w(\mathfrak{u},\mathfrak{v},\mathfrak{a})}[\mathfrak{bv}^{-1}],\quad\mathfrak{a},\mathfrak{u},\mathfrak{v}\in
\Psi_0\quad\mathrm{et}\quad\mathfrak{u}\mid\mathfrak{v}
\end{equation}
Il est clair que la restriction de la projection $[\bar \Psi_0]\longrightarrow \mathcal{A}$ \`a
$\bar \Psi_0$ est une
distribution ordinaire sur $K$. De plus toute distribution ordinaire 
$f\colon\bar
\Psi_0\longrightarrow A$ revient \`a se donner un
homomorphisme de groupes ab\'eliens 
$\widetilde f\colon\mathcal{A}\longrightarrow A$. Le groupe
$\mathcal{A}$ est
appel\'e la distribution ordinaire universelle sur $K$. Soit $\mathfrak{m}\in \Psi_0$ un id\'eal
entier, et soit
$\mathcal{A}_{\mathfrak m}$ l'image dans $\mathcal{A}$ du sous-groupe de $[\bar \Psi_0]$
engendr\'e par les symboles
$[\mathfrak{am}^{-1}],\mathfrak{a}\in \Psi_0$. Alors $\mathcal{A}_{\mathfrak m}$ est par
d\'efinition le groupe de
niveau $\mathfrak{m}$ de la distribution $\mathcal{A}$.
D'apr\`es l'article \cite{Y00}  le $\mathbb{Z}$-rang de
$\mathcal{A}_{\mathfrak m}$ est $\# G_{\mathfrak m}$,
ce qui peut se d\'eduire alternativement de l'\'etude qui suit.\par
Soit $\Delta$ le groupe ab\'elien libre engendr\'e par la r\'eunion disjointe

$$\Delta:=
 \langle \bigsqcup_{\mathfrak{n}\in \Psi_0}G_{\mathfrak n} \rangle ;$$
et consid\'erons son quotient $\hat{\mathcal{A}}:=\Delta/U$, o\`u $U$ est le groupe engendr\'e
par les sommes
formelles $S(\mathfrak{n},\mathfrak{p}^e,\sigma)$, $\mathfrak{n}\in \Psi_0$ et $\mathfrak{p}$
id\'eal maximal de
$\mathcal{O}_K$, d\'efinies ci-dessous~:
\begin{equation}\label{relred}
S(\mathfrak{n},\mathfrak{p}^e,\sigma):=\sigma-s(K_{\mathfrak{np}^e}/K_{\mathfrak{n}})\widetilde\sigma
\quad \mathrm{si}\quad\mathfrak{p}\mid\mathfrak{n}
\end{equation}
\begin{equation}\label{reluti}
S(\mathfrak{n},\mathfrak{p}^e,\sigma):=(1-(\mathfrak{p},K_{\mathfrak{n}}/K)^{-1})\sigma-s(K_{\mathfrak{np}^e}/K_{\mathfrak{n}})\widetilde\sigma
\quad \mathrm{si}\quad\mathfrak{p}\nmid\mathfrak{n}.
\end{equation}
Les automorphismes $\sigma\in G_{\mathfrak n}$ et $\widetilde\sigma\in G_{\mathfrak{np}^e}$ sont
choisis de telle
sorte que $\widetilde\sigma=\sigma$ sur $K_{\mathfrak{n}}$. Si
$\mathfrak{n}\mid\mathfrak{n}^\prime$ alors
$s(K_{\mathfrak{n}^\prime}/K_{\mathfrak n})$ est la somme dans $\Delta$
des \'el\'ements du groupe de Galois de
$K_{\mathfrak{n}^\prime}/K_{\mathfrak n}$, soit
$s(\mathrm{Gal}(K_{\mathfrak{n}^\prime}/K_{\mathfrak n}))$. Alors
on a
\begin{theorem}\label{iso}
Les groupes $\mathcal{A}$ et $\hat{\mathcal{A}}$ sont canoniquement isomorphes
en tant que $\mathbb{Z}[G]$-modules,
o\`u $G:=\limproj G_{\mathfrak n}$ est la limite projective des groupes $G_{\mathfrak n}$.
\end{theorem}
\noindent{\bf Preuve.} Remarquons que les relations (\ref{dist}) 
avec $\mathfrak{v}=\mathfrak{u}$ ou
$\mathfrak{v}=\mathfrak{up}$, $\mathfrak{p}$ \'etant un id\'eal premier,
engendrent toutes les autres. Si $\mathfrak{v}=\mathfrak{u}$ il s'agit de la relation
\begin{equation}\label{equiv}
[\mathfrak{au}^{-1}]-[\mathfrak{bv}^{-1}],\quad\mathfrak{b}\sim_{\mathfrak
u}\mathfrak{a}.
\end{equation}
Supposons que $\mathfrak{v}=\mathfrak{up}$ et soit $\mathfrak{a}\in\Psi_0$. Posons
$\mathfrak{d}:=pgcd(\mathfrak{u},\mathfrak{a})$, $\mathfrak{a}^\prime:=\mathfrak{ad}^{-1}$ et
$\mathfrak{n}:=\mathfrak{ud}^{-1}$. Consid\'erons alors l'entier $v_{\mathfrak
p}(\mathfrak{a}):=\mathrm{max}\{k\in\mathbb{N}\left\vert\ \mathfrak{p}^k\mid\mathfrak{a}\}
\right.$. Si $v_{\mathfrak p}(\mathfrak{a})<v_{\mathfrak p}(\mathfrak{u})$ alors la relation 
(\ref{dist}) s'\'ecrit simplement
\begin{equation}\label{dist1}
[\mathfrak{a}^\prime\mathfrak{n}^{-1}]-\sum_{\mathfrak{b}\in
w(\mathfrak{n},\mathfrak{np},\mathfrak{a}^\prime)}[\mathfrak{b}\mathfrak{n}^{-1}\mathfrak{p}^{-1}].
\end{equation}
Dans ce cas on a bien s\^ur $\mathfrak{p}\mid\mathfrak{n}$. En revanche si 
$v_{\mathfrak p}(\mathfrak{a})\geq v_{\mathfrak p}(\mathfrak{u})$ alors
$\mathfrak{p}\nmid\mathfrak{n}$. De plus dans
$\Psi^{-1}_{\mathfrak{n},\mathfrak{np}}(\psi_{\mathfrak n}(\mathfrak{a}^\prime))$ on trouve une
seule classe pouvant \^etre repr\'esent\'ee par un id\'eal non \'etranger \`a $\mathfrak{np}$.
Fixons en un et notons le $\mathfrak{b}_0$. Alors on obtient la relation
\begin{equation}\label{dist2}
[\mathfrak{a}^\prime\mathfrak{n}^{-1}]-[\mathfrak{b}_0\mathfrak{n}^{-1}\mathfrak{p}^{-1}]-\sum_{\overset {\scriptstyle \mathfrak{b}\in
w(\mathfrak{n},\mathfrak{np},\mathfrak{a}^\prime)}
{pgcd(\mathfrak{b},\mathfrak{np})=1}}[\mathfrak{b}\mathfrak{n}^{-1}\mathfrak{p}^{-1}].
\end{equation}
Ainsi l'application $[\mathfrak{an}^{-1}]\longrightarrow(\mathfrak{a}, K_{\mathfrak n}/K)$
 envoie (\ref{equiv}) sur $0$ tandis que l'image de (\ref{dist1}) et (\ref{dist2}) est la somme
$S(\mathfrak{n},\mathfrak{p},(\mathfrak{a}^\prime, K_{\mathfrak n}/K))$. On obtient donc un
homomorphisme de groupes ab\'eliens $\rho\colon\mathcal{A}\longrightarrow\hat{\mathcal{A}}$. Son
homomorphisme r\'eciproque $\rho^{-1}$ se d\'eduit de l'application
$\Delta\longrightarrow\mathcal{A}$ qui associe \`a $\sigma\in G_{\mathfrak n}$ la classe de
$[\mathfrak{an}^{-1}]$, o\`u $\mathfrak{a}\in\Psi_0$ est choisi de telle sorte que
$\sigma=(\mathfrak{a}, K_{\mathfrak n}/K)$. cette application est bien d\'efinie \'etant donn\'e
(\ref{equiv}).\qed\par
Notons $\Delta_{\mathfrak m}$ le sous groupe de $\Delta$ engendr\'e par $\lbrace \sigma\in
G_{\mathfrak
n},\mathfrak{n}\mid\mathfrak{m}\rbrace$, et
posons $U_{\mathfrak m}:= \Delta_{\mathfrak m}\cap U$. Alors $\mathcal{A}_{\mathfrak m}$ et
$\Delta_{\mathfrak m}/U_{\mathfrak m}$ sont isomorphes comme $\mathbb{Z}[G_{\mathfrak
m}]$-modules.  

\section{Sur la $\ell$-torsion de $\mathcal{A}_{\mathfrak m}$}\label{ltor}

 Il s'agit ici de trouver des conditions suffisantes pour que $\mathcal{A}_{\mathfrak m}$ soit
sans $\ell$-torsion.
Notre m\'ethode consiste \`a comparer la distribution universelle \`a la distribution d'Iwasawa
de $K$. Celle-ci
prend ses valeurs dans la limite directe $\Omega:=\limindu\mathbb{Q}[G_{\mathfrak n}]$ pour le
systeme inductif
$$\mathbb{Q}[G_{\mathfrak n}]\overset{\phi_{\mathfrak{n},\mathfrak{n}^\prime}}
{\hookrightarrow}
\mathbb{Q}[G_{\mathfrak{n}^\prime}],\quad\sigma\longmapsto s(K_{\mathfrak{n}^\prime}/K_{\mathfrak
n})\widetilde\sigma,$$
o\`u $\mathfrak{n}\mid\mathfrak{n}^\prime$ et $\widetilde\sigma$ est une extension quelconque de
$\sigma$ \`a
$K_{\mathfrak{n}^\prime}$. Pour introduire cette distribution, on 
d\'efinit d'abord les \'el\'ements d'Iwasawa suivants~:

\begin{equation}\label{alpha}
\alpha(\mathfrak{n},\mathfrak{n}^\prime):=s(\mathrm{Gal}(K_{\mathfrak{n}^\prime}/K_{\mathfrak
n}))\prod_{\mathfrak{p}\mid\mathfrak{n}}(1-\mathfrak p^*)
,\quad\mathfrak{n}\mid\mathfrak{n}^\prime,\end{equation}
o\`u $\mathfrak p^*:=\lambda_{\mathfrak p}^{-1}s(T_{\mathfrak p})/\# 
T_{\mathfrak p}$ est la moyenne des Frobenius inverse de
$\mathfrak p$, avec $T_{\mathfrak p}$:=$T_{\mathfrak
p}(\mathfrak{n}^\prime)$ et 
$\lambda_{\mathfrak p}\in
G_{\mathfrak{n}^\prime}/T_{\mathfrak p}$ est le Frobenius en $\mathfrak{p}$ dans
$K_{\mathfrak{n}^\prime}/K$.
Remarquons que l'on a
$$\phi_{\mathfrak{n}^\prime,\mathfrak{n}^{\prime\prime}}(\alpha(\mathfrak{n},\mathfrak{n}^\prime))=\alpha(\mathfrak{n},\mathfrak{n}^{\prime\prime}).$$
La distribution d'Iwasawa $\mathcal{F}\colon\Delta\longrightarrow\Omega$ peut \^etre d\'efinie comme
la limite des
applications $\mathcal{F}_{\mathfrak n}\colon\Delta_{\mathfrak
n}\longrightarrow\mathbb{Q}[G_{\mathfrak n}]$ o\`u pour
$\sigma\in G_{\mathfrak u}$, $\mathfrak{u}\mid\mathfrak{n}$ on a pos\'e
$$\mathcal{F}_{\mathfrak
n}(\sigma)=\phi_{\mathfrak{u},\mathfrak{n}}(\sigma\alpha(\mathfrak{u},\mathfrak{u}))=\widetilde\sigma\alpha(\mathfrak{u},\mathfrak{n}).$$
Consid\'erons maintenant le sous groupe $U(\mathfrak{n})$ de $\Delta_{\mathfrak n}$ engendr\'e
par les sommes
$S(\mathfrak{u},\mathfrak{p}^e,\sigma)$ pour lesquelles on a $\mathfrak{up}^e\mid\mathfrak{n}$.
On v\'erifie ais\'ement l'inclusion 
$U\subset\ker(\mathcal{F})$. On en d\'eduit
$U_{\mathfrak n}\subset\ker(\mathcal{F}_{\mathfrak n})$. D'autre part 
l'inclusion
$U(\mathfrak{n})\subset U_{\mathfrak n}$ donne une surjection
\begin{equation}\label{surj}
\ker(\mathcal{F}_{\mathfrak n})/U(\mathfrak{n})\twoheadrightarrow\ker(\mathcal{F}_{\mathfrak
n})/U_{\mathfrak
n}.
\end{equation}
On va voir avec (\ref{egn=m}) et le th\'eor\`eme \ref{zu} ci-dessous que 
ces deux groupes sont finis. Comme
$\mathrm{Im}(\mathcal{F}_{\mathfrak n})$ est
un $\mathbb{Z}$-module sans torsion, on pourra en d\'eduire les \'egalit\'es
$$\mathrm{Tor}(\Delta_{\mathfrak n}/U(\mathfrak{n}))=
\ker(\mathcal{F}_{\mathfrak n})/U(\mathfrak{n})\quad\mathrm{et}\quad
\mathrm{Tor}(\Delta_{\mathfrak n}/U_{\mathfrak
n})=\ker(\mathcal{F}_{\mathfrak n})/U_{\mathfrak n},$$
o\`u $\mathrm{Tor}(A)$ d\'esigne la torsion du groupe ab\'elien $A$.
En cons\'equence, si le quotient 
 $\Delta_{\mathfrak m}/U(\mathfrak{m})$ est sans 
$\ell$-torsion, il en est de m\^eme pour
$\Delta_{\mathfrak m}/U_{\mathfrak m}$. 
Cela incite \`a \'etudier le groupe
$\ker(\mathcal{F}_{\mathfrak
m})/U(\mathfrak{m})$, d'autant plus que la d\'efinition de 
$U(\mathfrak{m})$ est plus explicite que celle 
de $U_\mathfrak{m}$.\par
 Notons $\mathfrak{p}_1$, \ldots, $\mathfrak{p}_s$ les id\'eaux
premiers qui divisent
$\mathfrak{m}$, et $e_1$, \ldots, $e_s$ leurs exposants respectifs
 de sorte qu'on ait l'\'egalit\'e 
$\mathfrak{m}=\mathfrak{p}_1^{e_1}\cdots\mathfrak{p}_s^{e_s}$. Soit $\Sigma$ l'ensemble des
diviseurs
$\mathfrak{n}$ de $\mathfrak{m}$ tels que
$\mathfrak{n}=\mathfrak{p}_1^{t_1e_1}\cdots\mathfrak{p}_s^{t_se_s}$, avec $t_i\in\{0,1\}$.
On choisit sur $\Sigma$ une relation d'ordre totale $\prec$ prolongeant 
la division, c'est-\`a-dire  
v\'erifiant l'implication $(\mathfrak{n}\mid\mathfrak{n}^\prime)
 \Rightarrow ( \mathfrak{n}\prec\mathfrak{n}^\prime)$.
 La relation $\prec$ permet d'associer \`a tout
$\mathfrak{n}\in\Sigma$ les deux sous groupes de $\Delta_{\mathfrak m}$ 
d\'efinis par les syst\`emes de g\'en\'erateurs :
$$\left . \hat\Delta_{\mathfrak n}:=\langle\sigma\in G_{\mathfrak
u} \right \vert \quad\mathfrak{u}\in\Sigma\quad\mathrm{et}\quad\mathfrak{u}\prec\mathfrak{n}\rangle$$
$$\left . \hat{U}(\mathfrak{n}):=\langle S(\mathfrak{u},\mathfrak{p}_i^{e_i},\sigma)
\right \vert 
\quad\mathfrak{u}\in\Sigma,\quad\mathfrak{p}_i\nmid\mathfrak{u}\quad\mathrm{et}\quad\mathfrak{up}_i^{e_i}\prec\mathfrak{n}\rangle$$
La restriction de $\mathcal{F}_{\mathfrak m}$ \`a $\hat\Delta_{\mathfrak n}$ sera not\'ee
$\hat{\mathcal{F}_{\mathfrak n}}$. Pour $\mathfrak{n}=\mathfrak{m}$, les 
deux groupes  $\hat\Delta_{\mathfrak m}/\hat{U}(\mathfrak{m})$ et
$\Delta_{\mathfrak m}/U(\mathfrak{m})$ sont \'egaux, on a donc l'\'egalit\'e  
\begin{equation}\label{egn=m}\ker(\hat{\mathcal{F}_{\mathfrak m}})/\hat {U}({\mathfrak m})=
\ker(\mathcal{F}_{\mathfrak m})/U(\mathfrak{m}).\end{equation}
Nous utiliserons aussi l'id\'eal $S(\mathfrak{n})$ de
 $\mathbb{Z}[G_{\mathfrak n}]$
engendr\'e par les traces suivant les sous-groupes d'inertie, c'est-\`a-dire
les sommes 
$s(T_{\mathfrak{p}_i}(\mathfrak{n}))=s(K_{\mathfrak n}/K_{\mathfrak{n}_i})$,
 o\`u $\mathfrak{n}_i\mathfrak{p_i}^{e_i}=\mathfrak{n}$, et o\`u
 $\mathfrak{p}_i$ parcourt l'ensemble des id\'eaux premiers qui divisent
un $\mathfrak{n}\in\Sigma$ fix\'e. L'exposant du sous-groupe de torsion de $\mathbb{Z}[G_{\mathfrak
n}]/S(\mathfrak{n})$ sera
not\'e $z_{\mathfrak n}$.
\begin{theorem}\label{zu}
Soit $\mathfrak{n}\in\Sigma$, alors $\ker(\hat{\mathcal{F}_{\mathfrak n}})/\hat {U}({\mathfrak n})$
est fini et on a 
$$\Biggl(\prod_{\mathfrak{u}\prec\mathfrak{n}}z_{\mathfrak
u}\Biggr)\ker(\hat{\mathcal{F}_{\mathfrak
n}})\subset\hat{U}({\mathfrak n}).$$
\end{theorem}
\noindent{\bf Preuve.} On fait un raisonnement par r\'ecurrence sur $\mathfrak{n}$. Si
$\mathfrak{n}=(1)$ alors
$\ker(\hat{\mathcal{F}_{\mathfrak n}})=0$. Ainsi on peut supposer $\mathfrak{n}\neq(1)$. Soit alors
$\widetilde{\mathfrak{n}}$ le pr\'ed\'ecesseur de $\mathfrak{n}$, c'est-\`a-dire que
$\widetilde{\mathfrak{n}}=
max\{\mathfrak{u}\prec\mathfrak{n}, \mathfrak{u}\neq\mathfrak{n}\}$. On pose
$R:=\mathbb{Z}[G_{\mathfrak m}]$. La compos\'ee des deux applications naturelles
$$\hat\Delta_{\mathfrak n}\longrightarrow\mathbb{Z}[G_{\mathfrak n}] \quad\mathrm{et}\quad
\mathbb{Z}[G_{\mathfrak
n}]\longrightarrow\mathbb{Z}[G_{\mathfrak n}]/S(\mathfrak{n})$$ est un $R$-homomorphisme
surjectif avec un noyau
\'egal \`a $S(\mathfrak {n})+\hat\Delta_{\widetilde{{\mathfrak n}}}=\hat U({\mathfrak
n})+\hat\Delta_{\widetilde{{\mathfrak n}}}$.
Mais comme $\hat U(\widetilde{{\mathfrak n}})\subset\hat U(\mathfrak
{n})\cap\hat\Delta_{\widetilde{{\mathfrak
n}}}$ on obtient la suite exacte
$$\hat\Delta_{\widetilde{{\mathfrak n}}}/\hat U(\widetilde{{\mathfrak
n}})\longrightarrow\hat\Delta_{\mathfrak
n}/\hat U(\mathfrak {n})\longrightarrow\mathbb{Z}[G_{\mathfrak
n}]/S(\mathfrak{n})\longrightarrow0.$$
D'autre part $\mathcal{F}_{\mathfrak m}$ induit un homomorphisme de $R$-modules surjectif
$$\beta_{\mathfrak n}\colon\mathbb{Z}[G_{\mathfrak
n}]/S(\mathfrak{n})\longrightarrow\mathrm{Im}(\hat{\mathcal{F}_{\mathfrak
n}})/\mathrm{Im}(\hat{\mathcal{F}_{\widetilde{{\mathfrak n}}}}).$$
Ainsi, en applicant $\mathcal{F}_{\mathfrak m}$ \`a notre suite exacte on obtient le diagramme
commutatif~:
$$\begin{array}{ccccccccc}
 &   &\hat\Delta_{\widetilde{\mathfrak{n}}}/\hat U(\widetilde{\mathfrak {n}})&\hfl{}{}&
\hat\Delta_{\mathfrak n}/\hat U(\mathfrak{ n})&\hfl{}{}&\mathbb{Z}[G_{\mathfrak
n}]/S(\mathfrak{n})&\hfl{}{}&0\\
 &   &\vfl{}{}& &\vfl{}{}& &\vfl{}{\beta_{\mathfrak n}}& & \\
0&\hfl{}{}&\mathrm{Im}(\hat{\mathcal{F}_{\widetilde{\mathfrak{n}}}})&\hfl{}{}&\mathrm{Im}(\hat{\mathcal{F}_{\mathfrak
n}})&\hfl{}{}&\mathrm{Im}(\hat{\mathcal{F}_{\mathfrak n}})/\mathrm{Im}(\hat{
\mathcal{F}_{\widetilde{\mathfrak{n}}}})&\hfl{}{}&0
\end{array}$$
On en d\'eduit la suite exacte~:
\begin{equation}\label{kerbet}
\ker(\hat{\mathcal{F}_{\widetilde{\mathfrak{n}}}})/\hat {U}(\widetilde{\mathfrak
n})\longrightarrow\ker(\hat{\mathcal{F}_{\mathfrak n}})/\hat
{U}(\mathfrak{n})\longrightarrow\ker\beta_{\mathfrak
n}\longrightarrow 0.\end{equation}
Par r\'ecurrence $\ker(\hat{\mathcal{F}_{\widetilde{\mathfrak{n}}}})/\hat
{U}({\widetilde{\mathfrak{n}}})$ est fini et
est annul\'e par $\prod_{\mathfrak{u}\prec\widetilde{\mathfrak{n}}}z_{\mathfrak u}$. Il suffit
donc, pour conclure
la preuve du th\'eor\`eme, de v\'erifier que $\ker\beta_{\mathfrak n}$ est fini. Mais comme
$\ker\beta_{\mathfrak
n}$ est de type fini on doit seulement s'assurer que les modules $\mathbb{Z}[G_{\mathfrak
n}]/S(\mathfrak{n})$ et
$\mathrm{Im}(\hat{\mathcal{F}_{\mathfrak n}})/\mathrm{Im}(\hat{
\mathcal{F}_{\widetilde{{\mathfrak
n}}}})$ ont m\^eme
rang sur $\mathbb{Z}$. Or on peut prouver par la th\'eorie des caract\`eres que
$$\mathrm{rang}(S(\mathfrak{
n}))=\sum_{\overset {\scriptstyle \mathfrak{u}\in\Sigma,\mathfrak{u}\mid\mathfrak{n}}
{
\mathfrak{u}\neq\mathfrak{n}}}
(-1)^{d(\mathfrak{u},\mathfrak{n})+1}[K_{\mathfrak
u}:K],$$
o\`u $d(\mathfrak{u},\mathfrak{n}):=\#\{i,\ \mathfrak{p}_i\mid\mathfrak{n}\ \mathrm{et}\
\mathfrak{p}_i\nmid\mathfrak{u}\}$ 
(la preuve de \cite{JNT1}, Proposition 2.11 se retranscrit telle quelle).
 Il reste \`a calculer le rang de
$\mathrm{Im}(\hat{\mathcal{F}_{\mathfrak n}})/\mathrm{Im}(\hat{\mathcal{F}_{\widetilde{{\mathfrak {n}}}}})$.
Pour ce faire consid\'erons les modules suivants

\begin{definition}\label{XY} Soient $\mathfrak{u},\mathfrak{v}\in\Sigma$ tels que
$\mathfrak{u}\mid\mathfrak{v}$. On d\'efinit
alors les quatre $R$-modules $X_{\mathfrak u}, Y_{\mathfrak u}, X_{\mathfrak u}(\mathfrak{v})$ et
$Y_{\mathfrak
u}(\mathfrak{v})$ en donnant une famille de g\'en\'erateurs de chacun d'entre eux~:
$$\left . \overset {\quad} {X_{\mathfrak u}}
:=\langle\alpha({\mathfrak n},\mathfrak{m})\right \vert\
 \mathfrak{n}\in\Sigma,\
\mathfrak{n}\mid\mathfrak{u}\rangle,$$
$$\left . \overset {\quad} {Y_{\mathfrak u}}
:=\langle\alpha({\mathfrak n},\mathfrak{m})\right \vert\
\mathfrak{n}\in\Sigma,\
\mathfrak{n}\mid\mathfrak{u}\
\mathrm{et}\  \mathfrak{n}\neq\mathfrak{u}\rangle,$$
$$\left . \overset {\quad} {X_{\mathfrak u}(\mathfrak{v})}
:=\langle\alpha({\mathfrak n},\mathfrak{m})
\right \vert \ \mathfrak{n}\in\Sigma,\
\mathfrak{n}\mid\mathfrak{v}\ \mathrm{et}\ \mathfrak{n}\prec\mathfrak{u}\rangle,$$
$$\left . \overset {\quad} {Y_{\mathfrak u}(\mathfrak{v})}
:=\langle\alpha({\mathfrak n},\mathfrak{m})
\right \vert \ \mathfrak{n}\in\Sigma,\
\mathfrak{n}\mid\mathfrak{v},\ \mathfrak{n}\prec\mathfrak{u}\ \mathrm{et}\
\mathfrak{n}\neq\mathfrak{u}\rangle.$$
\end{definition}

\begin{proposition}\label{raindor} Soient $\mathfrak{u},\mathfrak{v}\in\Sigma$ tels que
$\mathfrak{u}\mid\mathfrak{v}$. Alors on a
$$\mathrm{rang}\Bigl(X_{\mathfrak u}(\mathfrak{v})/Y_{\mathfrak
u}(\mathfrak{v})\Bigr)=\mathrm{rang}\Bigl(X_{\mathfrak u}/Y_{\mathfrak u}\Bigr).$$
\end{proposition}
\noindent{\bf Preuve.} On voit que $X_{\mathfrak u}\subset X_{\mathfrak u}(\mathfrak{v})$ et
$Y_{\mathfrak
u}\subset Y_{\mathfrak u}(\mathfrak{v})$. De plus on a la suite exacte suivante 
$$0\longrightarrow X_{\mathfrak u}\cap Y_{\mathfrak u}(\mathfrak{v})/Y_{\mathfrak
u}\longrightarrow X_{\mathfrak
u}/Y_{\mathfrak u}\longrightarrow X_{\mathfrak u}(\mathfrak{v})/Y_{\mathfrak
u}(\mathfrak{v})\longrightarrow0.$$
Or le $\mathbb{Z}$-module $X_{\mathfrak u}\cap Y_{\mathfrak u}(\mathfrak{v})/Y_{\mathfrak u}$ est
fini puisqu'il
est de type fini et annul\'e par $[K_{\mathfrak v}:K_{\mathfrak u}]$. Ceci prouve la proposition
et montre en
particulier que $$\mathrm{rang}\Bigl(\mathrm{Im}(\hat{\mathcal{F}_{\mathfrak
n}})/\mathrm{Im}(\hat{\mathcal{F}_{\widetilde{\mathfrak
{n}}}})\Bigr)=\mathrm{rang}\Bigl(X_{\mathfrak n}/Y_{\mathfrak
n}\Bigr).$$
En effet on a $\mathrm{Im}(\hat{\mathcal{F}_{\mathfrak
n}})=X_{\mathfrak n}(\mathfrak{m})$ et $\mathrm{Im}(\hat{\mathcal{F}_{\widetilde{\mathfrak
{n}}}})=Y_{\mathfrak n}(\mathfrak{m})$.\qed

\noindent  Il nous faut maintenant calculer le rang du $\mathbb{Z}$-module $X_{\mathfrak
n}/Y_{\mathfrak n}$.
\begin{proposition}\label{radi}
On a $X_{\mathfrak n}=\phi_{\mathfrak{n},\mathfrak{m}}(\mathrm{Im}(\mathcal{F}_{\mathfrak n}))$.
De plus $\mathrm{Im}(\mathcal{F}_{\mathfrak n})$ est un $\mathbb{Z}$-module libre de rang $\# G_{\mathfrak n}$.
\end{proposition}
\noindent{\bf preuve.} La premi\`ere assertion est \'evidente. Pour montrer la seconde il suffit
de v\'erifier que si $\chi$ est un caract\`ere complexe du groupe $G_{\mathfrak n}$ alors
$\chi(\mathrm{Im}(\mathcal{F}_{\mathfrak n}))\neq0$. Or si on note $\mathfrak{n}_\chi$ le
conducteur de $\chi$ alors
$\chi(\alpha(\mathfrak{n}_\chi,\mathfrak{n}))=\#\mathrm{Gal}(K_{\mathfrak
n}/K_{\mathfrak{n}_\chi})$.\qed

Comme $\phi_{\mathfrak{n},\mathfrak{m}}$ est injective on a $\mathrm{rang}(X_{\mathfrak n})=\# G_{\mathfrak n}$.
 Par ailleurs on a $Y_{\mathfrak u}(\mathfrak{n})=X_{{\mathfrak
u}^\prime}(\mathfrak{n})$, o\`u
$\mathfrak{u}^\prime:=\mathrm{max}\{\mathfrak{t}\in\Sigma,\ \mathfrak{t}\mid\mathfrak{n},\
\mathfrak{t}\prec\mathfrak{u}\ \mathrm{et}\ \mathfrak{t}\neq\mathfrak{u}\}$. D'o\`u la relation
$$\sum_{\overset {\scriptstyle \mathfrak{u}\in\Sigma,\mathfrak{u}\mid\mathfrak{n}}
{\mathfrak{u}\neq\mathfrak{n}}}\mathrm{rang}\Bigl(X_{\mathfrak
u}(\mathfrak{n})/Y_{\mathfrak u}(\mathfrak{n})\Bigr)=\mathrm{rang}(Y_{\mathfrak n}),$$
qui peut \^etre utilis\'ee pour prouver, par r\'ecurrence sur $\vert
\mathfrak{n}\vert:=d((1),\mathfrak{n})$, la
formule
$$\mathrm{rang}\Bigl(X_{\mathfrak n}/Y_{\mathfrak
n}\Bigr)=\sum_{\mathfrak{u}\in\Sigma,\mathfrak{u}\mid\mathfrak{n}}(-1)^{d(\mathfrak{u},\mathfrak{n})}[K_{\mathfrak
u}:K].$$
Ainsi on a prouv\'e l'identit\'e
$$\mathrm{rang}\Bigl(\mathrm{Im}(\hat{\mathcal{F}_{\mathfrak n}})/\mathrm{Im}(\hat{\mathcal
{F}_{\widetilde{{\mathfrak
n}}}})\Bigr)=\mathrm{rang}\Bigl(\mathbb{Z}[G_{\mathfrak n}]/S({\mathfrak n})\Bigr),$$
et ceci compl\`ete la preuve de la proposition et du th\'eor\`eme.\qed
\begin{proposition}\label{deploie}
Posons $H:=K_{(1)}$ et soit $\ell$ un nombre premier. Soit $G_\ell(\mathfrak{m})$ le
$\ell$-Sylow du
groupe de Galois de $K_{\mathfrak m}/H$ et $G_{\mathfrak p}^{\ell}(\mathfrak{m})$, pour
$\mathfrak{p}\mid\mathfrak{m}$, le $\ell$-Sylow de $T_{\mathfrak p}(\mathfrak{m})$. Supposons
que
$G_\ell(\mathfrak{m})$ est le produit direct des groupes $G_{\mathfrak p}^{\ell}(\mathfrak{m})$.
Alors
$\ell\nmid\prod_{\mathfrak{u}\in\Sigma}z_{\mathfrak u}$. En particulier le groupe
$\mathcal{A}_{\mathfrak m}$ est
sans $\ell$-torsion.
\end{proposition}
\begin{rem} Il est possible de montrer que si $\mathfrak{u}$ est divisible par au plus deux
id\'eaux
premiers alors $z_{\mathfrak u}=1$.
\end{rem}
\noindent{\bf Preuve de la proposition \ref{deploie}.} Fixons $\mathfrak{u}\in\Sigma$. Les hypoth\`eses de la
proposition
impliquent que le $\ell$-Sylow $G_\ell(\mathfrak{u})$ de $\Gamma:=\mathrm{Gal}(K_{\mathfrak
u}/H)$ est le produit
direct des $\ell$-Sylow $G_{\mathfrak p}^{\ell}(\mathfrak{u})$ de $T_{\mathfrak
p}(\mathfrak{u}),\mathfrak{p}\mid\mathfrak{u}$. Cette remarque nous sera utile pour montrer que
$\ell\nmid
z_{\mathfrak u}$. La premi\`ere \'etape de la d\'emonstration consiste \`a d\'ecomposer
$\mathbb{Z}[G_{\mathfrak
u}]/S(\mathfrak{u})$ en somme directe de sous-$\mathbb{Z}[\Gamma]$-modules. En effet, les
$s(T_{\mathfrak
p}(\mathfrak{u})), \mathfrak{p}\mid\mathfrak{u}$ engendrent dans $\mathbb{Z}[\Gamma]$ un id\'eal
que nous noterons
$\widetilde S(\mathfrak{u})$. De plus si $\{\gamma_{\mathfrak u}(\tau),\
\tau\in\mathrm{Gal}(H/K)\}$ est un
syst\`eme de repr\'esentants de $G_{\mathfrak u}$ modulo $\Gamma$, alors on a l'isomorphisme de
$\mathbb{Z}[\Gamma]$-modules suivant~:
\begin{equation}\label{h=1}
\frac{\mathbb{Z}[G_{\mathfrak
u}]}{S(\mathfrak{u})}\simeq\bigoplus_{\tau\in\mathrm{Gal}(H/K)}\frac{\mathbb{Z}[\Gamma]}{\widetilde
S(\mathfrak{u})}\gamma_{\mathfrak u}(\tau).\end{equation}
Nous devons donc montrer que $\mathbb{Z}[\Gamma]/\widetilde S(\mathfrak{u})$ est sans
$\ell$-torsion. Cela
revient \`a v\'erifier que le $\mathbb{Z}_\ell[\Gamma]$-module
$$\mathbb{Z}_\ell\otimes_{\mathbb{Z}}\mathbb{Z}[\Gamma]/\widetilde
S(\mathfrak{u})\simeq\mathbb{Z}_\ell[\Gamma]/\widetilde S(\mathfrak{u})\mathbb{Z}_\ell[\Gamma]$$
est sans $\mathbb{Z}_\ell$-torsion. A cette fin on \'ecrit $\Gamma=G^\prime\times
G_\ell(\mathfrak{u})$, avec
$\ell\nmid\#G^\prime$, et on note $\mathfrak{X}$ l'ensemble des caract\`eres $\mathbb{Q}_\ell$-irr\'eductibles de
$G^\prime$. Soit $\{ e_\chi,\ \chi\in\mathfrak{X}\}$ le syst\`eme 
complet d'idempotents
 associ\'es aux \'el\'ements de $\mathfrak{X}$. Rappelons que ce syst\`eme 
v\'erifie les relations (orthogonalit\'e des caract\`eres) 
$\sum_{\chi\in\mathfrak{X}}e_\chi=1$ et $e_\chi e_{\chi^\prime}=0$ (resp. $e_\chi$) si
$\chi^\prime\neq\chi$ (resp. $\chi^\prime=\chi$). On en tire
l'isomorphisme suivant~:
$$\frac{\mathbb{Z}_\ell[\Gamma]}{\widetilde
S(\mathfrak{u})\mathbb{Z}_\ell[\Gamma]}\simeq\bigoplus_{\chi\in\mathfrak{X}}\frac{e_\chi\mathbb{Z}_\ell[\Gamma]}{e_\chi\widetilde
S(\mathfrak {u})\mathbb{Z}_\ell[\Gamma]}.$$
D'autre part, si pour un caract\`ere donn\'e $\chi\in\mathfrak{X}$ on pose
$A_\chi:=e_\chi\mathbb{Z}_\ell[G^\prime]$, alors on constate que les deux $A_\chi$-alg\`ebres
$$\frac{e_\chi\mathbb{Z}_\ell[\Gamma]}{e_\chi\widetilde S(\mathfrak {u})\mathbb{Z}_\ell[\Gamma]}\
\ \mathrm{et}\ \
A(\chi,\mathfrak{u}):=\frac{A_\chi[G_\ell(\mathfrak{u})]}{\widetilde
S(\mathfrak{u})A_\chi[G_\ell(\mathfrak{u})]}$$
sont isomorphes. De plus comme
$G_\ell(\mathfrak{u})=\mathop{\prod}\limits_{\mathfrak{p}\mid\mathfrak{u}}G_{\mathfrak
p}^{\ell}(\mathfrak{u})$ on
obtient la d\'ecomposi\-tion~:  
$$A(\chi,\mathfrak{u})\simeq\bigotimes_{\mathfrak{p}\mid\mathfrak{u}}\frac{A_\chi[G_{\mathfrak
p}^{\ell}(\mathfrak{u})]}{s(T_{\mathfrak p}(\mathfrak{u}))A_\chi[G_{\mathfrak
p}^{\ell}(\mathfrak{u})]}.$$
de $A(\chi,\mathfrak{u})$ en produit tensoriel sur $A_\chi$ des alg\`ebres
$$A(\chi,\mathfrak{u},\mathfrak{p}):=A_\chi[G_{\mathfrak p}^{\ell}(\mathfrak{u})]/s(T_{\mathfrak
p}(\mathfrak{u}))A_\chi[G_{\mathfrak p}^{\ell}(\mathfrak{u})].$$
Or les $A(\chi,\mathfrak{u},\mathfrak{p})$ sont des $A_\chi$-modules libres, et comme $A_\chi$
est un
$\mathbb{Z}_\ell$-module libre on d\'eduit que les $A(\chi,\mathfrak{u},\mathfrak{p})$ sont
eux-m\^emes
$\mathbb{Z}_\ell$-libres. Ce qui compl\`ete la preuve de la proposition.\qed
\begin{corollary}\label{Qli}
 Soit $\mathfrak{m}=n\mathbb{Z}$ un id\'eal de $\mathbb{Z}$. Alors le groupe de niveau
$\mathfrak{m}$ de
la distribution ordinaire universelle sur $\mathbb{Q}$ est un groupe ab\'elien libre.
\end{corollary}
\noindent{\bf Preuve.} Cela tient au fait que le groupe de Galois sur $\mathbb{Q}$ de l'extension
cyclotomique
$\mathbb{Q}(e^{\frac{2i\pi}{n}})$ est \'egal au produit direct des
groupes d'inertie.\qed \par
\begin{corollary}{inj}
Supposons que $\mathrm{Gal}(K_{\mathfrak m}/H)$ est \'egal au produit direct de ses
sous-groupes
d'inertie. Alors on a $U(\mathfrak{m})=U\cap\Delta_{\mathfrak m}$. De plus, pour tout diviseur
$\mathfrak{n}$ de
$\mathfrak{m}$, $\mathcal{A}_{\mathfrak n}$ est $\mathbb{Z}$-libre.
\end{corollary}
\section{Le cas des corps quadratiques imaginaires}\label{cohom}

Nous supposons \`a partir de maintenant et jusqu'\`a la fin de cet article
que $K$ est un corps quadratique imaginaire.
\subsection{Ennonc\'es des resultats.}
La proposition \ref{deploie}, appliqu\'ee 
\`a cette situation particuli\`ere, entra\^{\i}ne le r\'esul\-tat suivant
\begin{corollary}\label{lnmidw} Soit $w_K$ le nombre de racines de l'unit\'e du corps
 quadratique imaginaire $K$. Soit $\ell$ un nombre
premier qui ne divise pas $w_K$, $\ell\nmid w_K$. Alors la distribution ordinaire universelle sur
$K$ est sans
$\ell$-torsion.
\end{corollary}
\noindent{\bf Preuve.} Le lecteur peut facilement v\'erifier que le $\ell$-Sylow de $G_{\mathfrak
m}$ est \'egal au
produit direct des $\ell$-Sylow des groupes d'inertie $T_{\mathfrak p}(\mathfrak {m}),\
\mathfrak{p}\mid\mathfrak{m}$.\qed \par
Dans la section pr\'ec\'edente on a calcul\'e 
un multiple $\prod_{\mathfrak{u}\mid \mathfrak{m}} z_{\mathfrak u}$
de l'exposant de
$\mathrm{Tor}(\mathcal{A}_{\mathfrak m})$. Si on suit la preuve de la
proposition \ref{deploie}, on voit que les sous-groupes de torsion des quotients
$\mathbb{Z}[\mathrm{Gal}(K_{\mathfrak n}/H)]/\widetilde
S(\mathfrak{n})$, pour $\mathfrak{n}\mid\mathfrak{m}$, jouent un r\^ole 
pr\'epond\'erant. Pour le cas particulier qui nous int\'eresse 
en vue des applications de \cite{Ou98}, ces derniers 
sous-groupes de torsion sont d\'ecrits par le 

\begin{theorem}\label{toralg} Soit $\mathfrak{n}$ un \'el\'ement de $\Sigma$, 
 suppos\'e premier \`a $w_K$. Rappelons que 
$\vert\mathfrak{n}\vert$ d\'esigne
 le nombre d'id\'eaux premiers qui divisent $\mathfrak{n}$. On a
$$\mathrm{Tor}(\mathbb{Z}[\mathrm{Gal}(K_{\mathfrak n}/H)]/\widetilde
S(\mathfrak{n}))\simeq\left\lbrace
\begin{array}{cc}
0&\mathrm{si}\ \vert\mathfrak{n}\vert=1\ \mathrm{ou}\ \vert\mathfrak{n}\vert\ \mathrm{pair},\\
\mathbb{Z}/w_K\mathbb{Z}&\mathrm{sinon.}
\end{array}\right.$$
\end{theorem}
\begin{corollary}\label{borne}
Supposons que $\mathfrak{m}$ est un id\'eal propre de $\mathcal{O}_K$ et premier \`a
$w_K$. Alors
l'ordre du groupe $\mathrm{Tor}(\mathcal{A}_{\mathfrak m})$ divise $(w_K)^{ah}$,o\`u $h$ est le
nombre de classes
d'id\'eaux de $K$ et $a=2^{\vert\mathfrak{m}\vert-1}-\vert\mathfrak{m}\vert$.
\end{corollary}
\noindent{\bf Preuve du corollaire.} La suite exacte (\ref{kerbet}) permet de montrer par r\'ecur\-ren\-ce sur
$\mathfrak{n}\in\Sigma-\{(1)\}$
que l'ordre de $\ker(\hat{\mathcal{F}_{\mathfrak n}})/\hat {U}({\mathfrak n})$ divise
$(w_K)^{a_{\mathfrak n}h}$, o\`u 
$$a_{\mathfrak n}:=\#\{\mathfrak{u}\in\Sigma,\ \mathfrak{u}\prec\mathfrak{n},\
\vert\mathfrak{u}\vert\
\mathrm{impair}\ \mathrm{et}\ \vert\mathfrak{u}\vert\geq3\}.$$
En effet, $\ker\beta_{\mathfrak n}$ est un sous-groupe de $\mathrm{Tor}(\mathbb{Z}[G_{\mathfrak
n}]/S(\mathfrak{n}))$. D'o\`u, grace \`a (\ref{h=1}) et au th\'eor\`eme \ref{toralg},
 l'ordre de $\ker\beta_{\mathfrak n}$ divise
$(w_K)^h$. De plus $\ker\beta_{\mathfrak n}=\{0\}$, si $\vert\mathfrak{n}\vert$ est \'egal \`a
$1$ ou \`a un entier
pair. 
Il ne reste plus qu'\`a v\'erifier que $a_{\mathfrak m}=2^{\vert 
\mathfrak {m}\vert-1}-\vert \mathfrak{m}\vert$.\qed\par
Le reste de cette section est consacr\'e \`a la preuve du th\'eor\`eme \ref{toralg}.
\subsection {Preuve du th\'eor\`eme \ref{toralg}}
D'apr\`es le corollaire \ref{lnmidw}, il suffit de montrer 
que les $\ell$-parties de $\mathbb{Z}/w_K\mathbb{Z}$ et de 
$\mathrm{Tor}(\mathbb{Z}[\mathrm{Gal}
(K_{\mathfrak n}/H)]/\widetilde
S(\mathfrak{n}))$ co\"\i ncident 
pour les nombres premiers $\ell$ qui divisent $w_K$. Fixons donc un tel 
$\ell$ (en fait $\ell=2$ ou $3$).
Pour simplifier les notations on suppose que
$\mathfrak{n}=\mathfrak{m}$, ce qui ne
change rien \`a la g\'en\'eralit\'e 
de la d\'emonstration, et on pose $m=\vert\mathfrak m\vert$.
 On laisse de cot\'e le cas imm\'ediat 
$m=1$.
On commence par \'enoncer sous forme de remarques des faits
\'el\'ementaires qui nous serons utiles.
\begin{rem}\label{inergen}
Le groupe $\mathrm{Gal}(K_{\mathfrak m}/H)$ est engendr\'e par les sous-groupes d'inertie
$T_{\mathfrak p}(\mathfrak{m})$, $\mathfrak{p}\mid\mathfrak{m}$
(ce r\'esultat bien s\^ur vaut
pour tout corps de nombres).
\end{rem}
\begin{rem}\label{inermures} Comme $\mathfrak{m}$ est premier \`a $w_K$ et 
$m>1$ on a pour tout $i$ 
$$T_{\mathfrak{p}_i}(\mathfrak{m})\simeq(\mathcal{O}_K/
\mathfrak{p}_i^{e_i})^\times.$$
\end{rem}
\begin{rem}\label{inerdir} Soit $i_0\in\{1,\ldots,m\}$, alors les groupes d'inertie
$T_{\mathfrak{p}_i}(\mathfrak{m})$, $i\neq
i_0$ forment un produit direct.
\end{rem}
La remarque \ref{inermures} montre en particulier que le $\ell$-Sylow 
$G_{\mathfrak{p}_i}^\ell(\mathfrak{m})$
de
$T_{\mathfrak{p}_i}(\mathfrak{m})$ est cyclique. Posons alors $g_i:=
o(G_{\mathfrak{p}_i}^\ell(\mathfrak{m}))$ et
supposons que $g_m\leq g_i$, pour tout $i$. Les remarques \ref{inergen} et \ref{inerdir}
permettent de v\'erifier que l'on peut trouver
dans $\mathrm{Gal}(K_{\mathfrak m}/H)$ des \'el\'ements $\tau_1,\ldots,\tau_m$
tels que
$G_{\mathfrak{p}_i}^\ell(\mathfrak{m})=\langle\tau_i\rangle,$ pour tout
$i\in\{1,\ldots,m-1\}$, et
$G_\ell(\mathfrak{m})=\langle\tau_1\rangle\times\langle\tau_2\rangle\times
\cdots\times\langle\tau_m\rangle$. De plus on a
$G_{\mathfrak{p}_m}^\ell(\mathfrak{m})=\langle j\rangle$, o\`u 
$$j:= \prod_{i=1}^m\tau_i^{\frac{g_i}{g_m}}$$
Notons que l'ordre de $\tau_m$ est $o(\langle\tau_m\rangle)=g_m/\ell^r$, o\`u $\ell^r$ est la plus grande puissance de $\ell$ divisant $w_K$. Ecrivons, comme en section \ref{ltor},
$\mathrm{Gal}(K_{\mathfrak m}/H)=G^\prime\times
G_\ell(\mathfrak{m})$ et soit
$G^\prime_i$ le sous-groupe de $G^\prime$ tel que
$T_{\mathfrak{p}_i}(\mathfrak{m})=G^\prime_i\times
G_{\mathfrak{p}_i}^\ell(\mathfrak{m})$.\par
Il s'agit de calculer la $\ell$-torsion de $A(\chi,\mathfrak{m})$ 
pour tout caract\`ere $\mathbb{Q}_\ell$-irr\'educti\-ble du groupe
$G^\prime$. On a d'abord le lemme \'evident :
\begin{lemma}\label{obv} Soit $\Lambda:=\mathbb{Z}[G_\ell(\mathfrak{m})]$ et $\Lambda_\chi$ le
$\Lambda$-module engendr\'e par
les traces
$s(G_{\mathfrak{p}_i}^\ell(\mathfrak{m}))$ pour lesquelles $\chi$ est trivial sur
$G^\prime_i$. Alors
$A(\chi,\mathfrak{m})$ est naturellement isomorphe \`a
$A_\chi\otimes_{\mathbb{Z}}\Lambda/\Lambda_\chi$.
\end{lemma}
\begin{proposition}\label{lib} Soit $P$ une partie de $\{1,\ldots,m\}$ distincte de
$\{1,\ldots,m\}$. 
Notons $\Lambda(P)$ le $\Lambda$-module
engendr\'e par les
traces $s(\langle\tau_i\rangle)$, $i\in P$. Alors le module $\Lambda/\Lambda(P)$ est libre sur
$D:=\mathbb{Z}[\langle\tau_k,\ k\not\in P\rangle]$. En particulier si $\chi$ n'est pas trivial sur
$G^\prime_m$
 alors $\Lambda/\Lambda_\chi$ est sans $\mathbb{Z}$-torsion.
\end{proposition}
\noindent{\bf Preuve.} Il est clair que $\Lambda/\Lambda(P)$ est isomorphe au
 produit tensoriel sur $D$
$$\bigotimes_{n\in P}\frac{D[\langle\tau_n\rangle]}{s(\langle\tau_n\rangle)D}.$$
Cela donne la premi\`ere assertion puisque 
$D[\langle\tau_n\rangle]/s(\langle\tau_n\rangle)D$ est libre sur $D$ (de rang $g_n-1$).\par
D'autre part pour $i\neq m$ les $\tau_i$ engendrent 
les $\ell$-groupes d'inerties : $\langle \tau_i\rangle = 
G_{\mathfrak p_i}^\ell(\mathfrak m )$. En particulier 
lorsque $\chi$ est non trivial sur $G_m^\prime$, on a l'\'egalit\'e 
$\Lambda_\chi=\Lambda(P_\chi)$ pour 
$\left . \overset{\quad} {P_\chi}=\{i\right \vert \chi(G^\prime_i)=\{1\}\}$.
Et la seconde assertion en d\'ecoule.\qed \par
Pour \'etudier la torsion de 
$\mathbb{Z}[\mathrm{Gal}(K_{\mathfrak n}/H)]/\widetilde
S(\mathfrak{n})$, on est amen\'e \`a consid\'erer 
les $\Lambda$-modules $\Theta(P)$ 
engendr\'e par les $s(G_{\mathfrak{p}_i}^\ell(\mathfrak{m}))$,
$i\in P$. Remarquons que si $m\notin P$, alors $\Theta(P)$ est \'egal 
au module $\Lambda(P)$ d\'efini ci-dessus.
En particulier, toujours si $m\notin P$, on a l'isomorphisme
$$\Lambda/\Theta(P\cup\{s\})\simeq \frac { \Lambda/\Lambda(P)}
{s(\langle j\rangle)\Lambda/\Lambda(P)}.$$
De plus comme $\Lambda/\Lambda(P)$ est $\mathbb{Z}$-libre d'apr\`es la 
proposition \ref{lib} ci-dessus on a
\begin{equation}\label{tor=h2}
\mathrm{Tor}(\Lambda/\Theta(P\cup\{s\}))\simeq H^2(\langle j\rangle,
\Lambda/\Lambda(P)).\end{equation}
La preuve du th\'eor\`eme \ref{toralg} consiste \`a calculer ces groupes de cohomologie.
\begin{proposition}\label{hpq=0} Soit $P\subset Q$ et $P\neq Q$ deux parties de 
$\{1,\ldots,m\}$ telles que $m\not\in Q$. Soit $D_Q$ le
sous-groupe de $G_\ell(\mathfrak{m})$ engendr\'e par les $G_{\mathfrak{p}_i}^\ell(\mathfrak{m})$,
$i\not\in Q$.
Alors on a
\begin{equation}\label{Hpq}
H^k_{P,Q}:=H^k(D_Q,\Lambda/\Lambda(P))=0,\ \mathrm{pour\ tout}\ k\geq1.\end{equation}
\end{proposition}
\noindent{\bf Preuve.} Nous allons d\'emontrer la proposition par r\'ecurrence sur le cardinal de
$P$. Le cas o\`u
$P=\emptyset$ est trivial puisqu'alors $\Lambda/\Lambda(P)=\Lambda$
qui est libre sur tout sous-groupe de
$G_\ell(\mathfrak{m})$. Supposons que (\ref{Hpq}) est prouv\'ee pour tout
couple $(P^\prime,Q^\prime)$ de sous-ensembles
de $\{1,\ldots,m\}$ v\'erifiant les hypoth\`eses de la proposition et tels que
le cardinal de $P^\prime< \#P$. On
peut bien s\^ur \'ecrire $P$ comme une r\'eunion disjointe $P=\{i_0\}\cup P^\prime$ o\`u $i_0\in
P$. On a alors la suite exacte
$$0\longrightarrow
s(\langle\tau_{i_0}\rangle)(\Lambda/\Lambda(P^\prime))\longrightarrow\Lambda/\Lambda(P^\prime)\longrightarrow\Lambda/\Lambda(P)\longrightarrow0.$$
D'o\`u l'on tire
$$H^k_{P^\prime,Q}\longrightarrow H^k_{P,Q}\longrightarrow
H^{k+1}(D_Q,s(\langle\tau_{i_0}\rangle)(\Lambda/\Lambda(P^\prime)))
\longrightarrow
H^{k+1}_{P^\prime,Q}.$$
L'hypoth\`ese de r\'ecurrence entra\^{\i}ne que les deux groupes $H^k_{P^\prime,Q}$ et
$H^{k+1}_{P^\prime,Q}$ sont
nuls. Il vient donc
$$
H^k_{P,Q}\simeq H^{k+1}(D_Q, s(\langle\tau_{i_0}\rangle)(\Lambda/
\Lambda(P^\prime))).$$
Or $\Lambda/\Lambda(P^\prime)$
est libre sur $\langle\tau_{i_0}\rangle$ comme \'enonc\'e dans la
proposition \ref{lib}. En
particulier on a $s(\langle\tau_{i_0}\rangle)(\Lambda/\Lambda(P^\prime))
=(\Lambda/\Lambda(P^\prime))^{\langle\tau_{i_0}\rangle}$ et 
$$H^n(\langle\tau_{i_0}\rangle,\Lambda/\Lambda(P^\prime))=
H^n(\langle\tau_{i_0}\rangle,\Lambda(P^\prime))=0\ \mathrm{pour\ tout}\ n>0.$$
D'autre part si on pose $Q^\prime:=Q-\{i_0\}$ alors on a $D_{Q^\prime}=\langle 
D_Q,\tau_{i_0}\rangle$. Comme on a
suppos\'e $P\neq Q$, on a $\{i_0\}\subsetneq Q$.
La remarque \ref{inergen} montre donc l'\'egalit\'e
 $D_Q\cap\langle\tau_{i_0}\rangle=\{1\}$, qui permet d'identifier
$D_Q$ avec le quotient $D_{Q^\prime}/\langle\tau_{i_0}\rangle$.
D'apr\`es le corollaire p. 118 de \cite{Self} ou le Theorem 2 p. 161 de \cite{Lang67}, l'inflation en degr\'e
$k+1$ donne donc un isomorphisme : 
$$H^{k+1}(D_Q,(\Lambda/\Lambda(P^\prime))^{\langle\tau_{i_0}\rangle})\simeq
 H^{k+1}(D_{Q^\prime},\Lambda/\Lambda(P^\prime)).$$
Mais comme $\# P^\prime < \# P$, 
le terme de droite est trivial, d'o\`u la proposition.\qed
\begin{corollary}\label{chitriv} Si le caract\`ere $\chi$ est non-trivial sur 
$G^\prime$
alors le $\mathbb{Z}$-module $\Lambda/\Lambda_\chi$ est sans
torsion.
\end{corollary}
\noindent{\bf Preuve}. Si $\chi$ est non trivial sur $G^\prime_m$,
c'est la proposition \ref{lib}. Sinon soit $P$ l'ensemble des indices $i\neq m$ tels
que $\chi$ est trivial sur
$G_i^\prime$ et soit $Q:=\{1,\ldots,m-1\}$. Comme $\chi$ est suppos\'e 
trivial sur $G^\prime_m$ et non sur $G^\prime$ on a l'inclusion stricte
$P\subsetneq Q$. D'apr\`es la d\'efinition de $\Theta$ et (\ref{tor=h2}) on a 
$$\mathrm{Tor}(\Lambda/\Lambda_\chi)=\mathrm{Tor}(\Lambda/\Theta(P\cup\{s\}))
\simeq H^2(D_Q,\Lambda/\Lambda(P)).$$
La proposition \ref{hpq=0} permet de conclure.\qed\par
Il reste \`a calculer la torsion de 
$\Lambda/\Lambda_\chi$ lorsque $\chi$ est trivial. 
Dans ce cas on a $\Lambda_\chi=\Theta(P_m)$ o\`u 
$P_i:=\{1,\ldots,i\}$. D'apr\`es (\ref{tor=h2}) on doit calculer le
groupe $H^2(\langle j\rangle,\Lambda/\Lambda(P_{m-1}))$.
La d\'emarche suivie pour prouver la proposition \ref{hpq=0}
montre en fait que l'on a
$$ H^k(D_{P_{i+1}},\Lambda/\Lambda(P_{i+1}))\simeq H^{k+1}(D_{P_i},
\Lambda/\Lambda(P_i)),\ i\in\{1,\ldots,m-2\}.$$
En particulier on a
$$H^2(\langle j\rangle,\Lambda/\Lambda(P_{m-1}))
\simeq H^m(D_{P_1},\Lambda/\Lambda(P_1))\simeq
H^{m+1}(D_{P_1},\Lambda^{\langle\tau_1\rangle}).$$
Or $\Lambda^{\langle\tau_1\rangle}=s(\langle\tau_1\rangle)\Lambda$
est libre sur le sous-groupe $\langle\tau_2,\ldots,\tau_{m-1}\rangle$ de $D_{P_1}$. L'inflation en degr\'e $m+1$ donne alors l'isomorphisme 
$$H^{m+1}(D_{P_1},\Lambda^{\langle\tau_1\rangle})\simeq H^{m+1}(\langle j\rangle,\Lambda^{\langle\tau_1,\ldots,\tau_{m-1}\rangle}).$$
Le groupe de cohomologie de droite se calcul sans difficulte: il est \'egal \`a $0$ si $m$ est pair, et a $\mathbb{Z}/\ell^r\mathbb{Z}$ si $m$ est impair. Ceci ach\`eve la preuve du th\'eor\`eme \ref{toralg}.

\section{Construction d'exemples o\`u $\mathrm{Tor}(\mathcal{A}_{\mathfrak
m})\neq\{0\}$}\label{contrex}
Dans cette section $K$ est toujours un corps quadratique imaginaire.
Nous avons vu pr\'ec\'edemment que $\mathcal{A}_{\mathfrak m}$ est sans
torsion si $s\leq2$ et que pour $s=3$, 
$\#\mathrm{Tor}(\mathcal{A}_{\mathfrak m})$ divise $(w_K)^h$.

\begin{proposition}\label{torsex} On suppose que $w_K=2$ 
et que $K$ contient trois id\'eaux premiers principaux $\mathfrak p_1 ,
\mathfrak p_2, \mathfrak p_3$ tels que $ N \mathfrak p_i\equiv 3 \ {\rm modulo }\ 4$.
On pose $\mathfrak m=\mathfrak p_1\mathfrak p_2\mathfrak p_3 $. Alors $\mathrm{Tor}
(\Delta_{\mathfrak m}/U(\mathfrak m))\cong (\mathbb {Z}/2\mathbb{Z})^h$ et $\mathrm {Tor}
(\mathcal{A}_{\mathfrak m})\neq0$.
\end{proposition}
Pour r\'ealiser les conditions de cette proposition, il suffit 
par exemple de choisir $K$ tel que $\sqrt{-1}\not\in H$ (ce qui se produit 
si $2$ n'est pas ramifi\'e dans $K$)
et de choisir trois id\'eaux premiers $\mathfrak {p}_i$ tels que
$(\mathfrak{p}_i,H(\sqrt{-1})/K)$ soit \'egal \`a l'\'el\'ement non trivial de
$\mathrm{Gal}(H(\sqrt{-1})/H)$.
\par
\noindent {\bf Preuve de la proposition \ref{torsex}.} 
Lorsque $m=3$ on a n\'ecessairement :

\noindent  $\mathrm{Tor}(\Delta_{\mathfrak m} / U(\mathfrak
m))\simeq\ker(\beta_{\mathfrak m})$. En effet il suffit de r\'e\'ecrire la suite exacte 
(\ref{kerbet}) en
remarquant que les probl\`emes de torsion ne se manifestent pas avant $\mathfrak{m}$. Mais vu la
d\'ecomposition (\ref{h=1}) et le th\'eor\`eme \ref{toralg}
il suffit simplement de v\'erifier  les deux
conditions $\ker(\beta_{\mathfrak m})=\mathrm{Tor} (\mathbb{Z}[G_{\mathfrak m}]/S(\mathfrak{m}))$
et $\mathrm{Tor}(\mathcal{A}_{\mathfrak m})\neq0$. C'est le cas si le seul \'el\'ement de torsion
de $\mathbb{Z}[\mathrm{Gal}(K_{\mathfrak m}/H)]/\widetilde S(\mathfrak{m})$ appartient \`a
$\ker(\beta_{\mathfrak m})$ d'une part, et que son ant\'ec\'edent par l'isomorphisme 
(\ref{kerbet})
n'appartient pas \`a $U_{\mathfrak m}/U(\mathfrak{m})$ d'autre part.
On a bien s\^ur $\ell=2$,
$o(G_{\mathfrak{p}_i}^\ell(\mathfrak{m}))=2$ et $o(G_\ell(\mathfrak{m}))=4$. On peut donc poser
$G_{\mathfrak{p}_i}^\ell(\mathfrak{m})=\langle\tau_i\rangle, i=1,2,3$, avec
$\tau_3=\tau_1 \tau_2$ et
$G_\ell(\mathfrak{m})=\langle\tau_1\rangle\times\langle\tau_2\rangle$.\par
Si $\chi$ est le caract\`ere trivial de $G^\prime$ alors $\Lambda/\Lambda_\chi$ est isomorphe \`a
$\mathbb{Z}/2\mathbb{Z}$, d'apr\`es le th\'eor\`eme \ref{toralg}.
Comme $A_\chi=s(G^\prime)\mathbb{Z}_\ell[G^\prime]$
est libre de rang 1 sur $\mathbb{Z}_\ell$, le seul \'el\'ement de
torsion de $A(\chi,\mathfrak{m})$ est la classe de $s(G^\prime)$.
Dire que cette classe appartient \`a $\ker\beta_{\mathfrak m}$ revient simplement \`a v\'erifier
que $s(G^\prime) \alpha(\mathfrak m,\mathfrak m)$ s'exprime en fonction des
$\alpha(\mathfrak{n},\mathfrak{m})$ pour
$\mathfrak{n}\mid\mathfrak{m}$ et $\mathfrak{n}\neq\mathfrak{m}$. On a 
$$ 2s(G^\prime) \alpha(\mathfrak m,\mathfrak m)=
((1+\tau_1)+(1+\tau_2)-\tau_1(1+\tau_3))s(G^\prime) 
\alpha(\mathfrak m,\mathfrak m).
$$
On va v\'erifier que chaque terme de la forme 
$(1+\tau_i) s(G^\prime) \alpha (\mathfrak m, \mathfrak m)$ appartient en fait 
\`a $2\mathbb{Z}[\mathrm {Gal}(K_{\mathfrak m}/H)] \ \alpha(\mathfrak {m}/
\mathfrak{p}_i, \mathfrak{m})$, en simplifiant par $2$ cela donne 
 l'expression attendue pour $s(G^\prime) \alpha(\mathfrak m, \mathfrak m)$.
Comme en (\ref{alpha}), on note $\lambda_{\mathfrak p_i}$ un relev\'e dans
$G_{\mathfrak m}$ du Frobenius en $\mathfrak p_i$ et on pose  
$\mathfrak p_i^*=\lambda^{-1}_{\mathfrak p_i}s(T_{\mathfrak p_i})/\# 
T_{\mathfrak p_i}$.
 Rappelons aussi que $G^\prime=G^\prime_1\times G^\prime_2\times G^\prime_3$
 o\`u les $G^\prime_i$ sont d\'efinis par $T_{\mathfrak p_i}=G^\prime_i
\times\{1,\tau_i\}$. D'o\`u l'\'egalit\'e  
$$(1+\tau_i)s(G^\prime)\alpha(\mathfrak m,\mathfrak m)
=s(T_{\mathfrak{p}_i}\times G^\prime_j\times G^\prime_k)(1-\lambda_{\mathfrak p_i}^{-1}) 
(1-\mathfrak p_j^*)(1-\mathfrak p_k^*)$$
Par hypoth\`ese l'automorphisme $\lambda_{\mathfrak p_i}\in\Gamma:=\mathrm{Gal}(K_{\mathfrak
m}/H)$. De plus le groupe $\Phi:=T_{\mathfrak{p}_i}\times G^\prime_j\times G^\prime_k$ \'etant
d'indice 2 dans $\Gamma$ on obtient
$$s(\Phi)(1-\lambda_{\mathfrak p_i}^{-1})=\left\lbrace
\begin{array}{cc}
0&\mathrm{si}\ \lambda_{\mathfrak p_i}\in\Phi\\
s(\Gamma)-2\lambda_{\mathfrak p_i}^{-1}s(\Phi)&\mathrm{sinon.}
\end{array}\right.$$ 
Or $s(\Gamma)(1-\mathfrak p_j^*)=s(\Gamma)(1-\lambda_{\mathfrak p_j}^{-1})=0$ 
puisque
$\lambda_{\mathfrak p_j}\in\Gamma$. Il suffit maintenant de remarquer que 
$s(\Phi)(1-\mathfrak p_j^*)(1-\mathfrak p_k^*)=s(G^\prime_j\times
G^\prime_k)\alpha(\mathfrak{p}_j\mathfrak{p}_k,\mathfrak{m})$ pour pouvoir
conclure que $s(G^\prime) \alpha(\mathfrak m,\mathfrak m)$ s'exprime
bien en fonction des
$\alpha(\mathfrak{n},\mathfrak{m})$ pour
$\mathfrak{n}\mid\mathfrak{m}$ et $\mathfrak{n}\neq\mathfrak{m}$. Plus pr\'ecis\'ement on a mis
en \'evidence des sommes $x_i\in\mathbb{Z}[\mathrm{Gal}(K_{\mathfrak{p}_j\mathfrak{p}_k}/K)]\subset \Delta_{\mathfrak m/\mathfrak{p}_i}$
telles que
$$R:=s(G^\prime)+x_1+x_2+x_3\in\ker(\mathcal{F}).$$
Pour d\'emontrer la proposition, il nous reste \`a prouver que cet \'el\'ement n'appartient pas \`a $U$. Cela peut se faire au
moyen d'un crit\`ere de parit\'e. Pour ce, on
introduit l'homomorphisme $\nu :\Delta \longrightarrow\mathbb{Z} $
d\'efini par $\nu(\sigma)=1$ si $\sigma\in G_{\mathfrak n}$ pour un $\mathfrak n$ 
de la forme $\mathfrak n = \prod_{i=1}^3 \mathfrak{p}_i^{e_i}$ avec des $e_i\geq 1$, et 
$\nu(\sigma)=0$ sinon. Clairement $\nu(R)=\# G^\prime $ est
impair tandis que 
\begin{lemma}\label{parite} $\nu(U)\subset 2\ \mathbb{Z}$.\end{lemma}
\noindent{\bf Preuve du lemme.} On \'etudie les images par $\nu$ des
g\'en\'erateurs de U not\'es $S(\mathfrak{n},\mathfrak{p}^e,\sigma)$ en
fin de section \ref{thmgen}.
Pour les g\'en\'erateurs du type (\ref{reluti}) la parit\'e est imm\'ediate
puisque $[K(\mathfrak{np}):K(\mathfrak{n})]$ est pair d\`es que 
$\mathfrak{p} \nmid \mathfrak{n}$ et $\mathfrak{n}\neq (1)$.
En ce qui concerne les g\'en\'erateurs du type (\ref{relred}) leur image
par $\nu$ est nulle sauf lorsque $\sigma \in G_{\mathfrak n}$ pour les $\mathfrak n$ de la 
forme $\mathfrak {n}=\prod_{i=1}^3 \mathfrak {p}_i^{e_i}$ avec des $e_i\geq 1$. Dans ce dernier cas
le $\mathfrak p$ de (\ref{relred}) est l'un des trois $\mathfrak p_i$,
et $[K(\mathfrak {np}^e):
K(\mathfrak {n})]=N(\mathfrak {p^e})$ est impair. 
D'o\`u la parit\'e de $\nu(S(\mathfrak{n},\mathfrak{p}^e,\sigma))=
1+[K(\mathfrak {np}^e):
K(\mathfrak {n})]$.\qed




\begin{thebibliography}{1}

\bibitem{JNT1}
Jean-Robert Belliard.
\newblock Sur la structure galoisienne des unit\'es circulaires dans les
  ${{\mathbb Z}}\sb p$-extensions.
\newblock {\em J. Number Theory}, 69(1):16--49, 1998.

\bibitem{DR80}
Pierre Deligne and Kenneth~A. Ribet.
\newblock Values of abelian ${L}$-functions at negative integers over totally
  real fields.
\newblock {\em Invent. Math.}, 59(3):227--286, 1980.

\bibitem{Ku79}
Daniel~S. Kubert.
\newblock The universal ordinary distribution.
\newblock {\em Bull. Soc. Math. France}, 107(2):179--202, 1979.

\bibitem{Lang67}
Serge Lang.
\newblock {\em Rapport sur la cohomologie des groupes}.
\newblock W. A. Benjamin, Inc., New York-Amsterdam, 1967.

\bibitem{MSD}
B.~Mazur and P.~Swinnerton-Dyer.
\newblock Arithmetic of {W}eil curves.
\newblock {\em Invent. Math.}, 25:1--61, 1974.

\bibitem{Ou98}
Hassan Oukhaba.
\newblock On elliptic units of abelian extensions of a given imaginary
  quadratic field.
\newblock {\em Pr\'epublication}, 1998.

\bibitem{Self}
Jean-Pierre Serre.
\newblock {\em Local fields}.
\newblock Springer-Verlag, New York, 1979.

\bibitem{Wa}
Lawrence~C. Washington.
\newblock {\em Introduction to cyclotomic fields}.
\newblock Springer-Verlag, New York, second edition, 1997.

\bibitem{Y00}
Linsheng Yin.
\newblock Distributions on a global field.
\newblock {\em J. Number Theory}, 80(1):154--167, 2000.

\end{thebibliography}

\end{document}